\documentclass[12pt,reqno]{amsart}
\pdfoutput=1

\usepackage{amssymb}
\usepackage{amsthm}
\usepackage{amsmath}
\usepackage{latexsym}
\usepackage[T1]{fontenc}
\usepackage[utf8]{inputenc}
\usepackage[english]{babel}

\usepackage{graphicx}
\usepackage{wrapfig}
\usepackage[justification=centering, labelfont=bf]{caption}
\usepackage{mathtools}
\usepackage{amsbsy}
\usepackage[inline]{enumitem}
\usepackage{mathrsfs}
\usepackage{array}
\usepackage{multicol}
\usepackage{stmaryrd}
\usepackage{cancel}
\usepackage{lmodern}
\usepackage{mathabx}
\usepackage{upgreek}
\usepackage{titlesec}
\usepackage{titletoc}
\usepackage[spacing=true,kerning=true,babel=true,tracking=true]{microtype}
\usepackage[foot]{amsaddr}
\usepackage[left=1in,right=1in,top=1in,bottom=1in,bindingoffset=0cm]{geometry}
\usepackage{bm}
\usepackage{centernot}
\usepackage{mdframed}
\usepackage{xspace}
\usepackage[skins]{tcolorbox}
\usepackage{framed}
\usepackage{float}
\usepackage{textgreek}

\usepackage[unicode=true,
 bookmarks=false,
 breaklinks=true,pdfborder={0 0 0},backref=false,colorlinks=false]
 {hyperref}
 
\date{}

\newtheoremstyle{bfnote}%
{}{}%
{\slshape}{}%
{\bfseries}{\bfseries.}%
{ }%
{\thmname{#1}\thmnumber{ #2}\thmnote{ \ep{\normalfont{}#3}}}

\theoremstyle{bfnote}
\newtheorem{theo}{Theorem}[section]
\newtheorem*{theo*}{Theorem}

\newtheorem{lemma}[theo]{Lemma}

\newtheorem{corl}[theo]{Corollary}
\newtheorem{conj}[theo]{Conjecture}

\newtheorem*{corl*}{Corollary}

\theoremstyle{definition}
\newtheorem{defn}[theo]{Definition}
\newtheorem*{defn*}{Definition}

\newtheorem*{exmp*}{Example}
\newtheorem{prob}[theo]{Problem}

\theoremstyle{remark}
\newtheorem*{ques*}{Question}
\newtheorem*{remk*}{Remark}

\newcommand*{\myproofname}{Proof}

\makeatletter
\newcommand{\neutralize}[1]{\expandafter\let\csname c@#1\endcsname\count@}
\makeatother

\newcommand{\R}{\mathbb{R}}

\renewcommand{\P}{\mathbb{P}}
\renewcommand{\epsilon}{\varepsilon}
\newcommand{\eps}{\epsilon}
\renewcommand{\phi}{\varphi}
\renewcommand{\theta}{\vartheta}
\renewcommand{\leq}{\leqslant}
\renewcommand{\geq}{\geqslant}
\renewcommand{\le}{\leq}
\renewcommand{\ge}{\geq}

\newcommand{\bemph}[1]{{\normalfont#1}} %
\newcommand{\ep}[1]{\bemph{(}#1\bemph{)}} %

\numberwithin{equation}{section}

\titleformat{\section}[block]{\large\bfseries\sffamily}{\thesection.}{1ex}{}
\titleformat{\subsection}[block]{\bfseries\sffamily}{\thesubsection.}{1ex}{}
\titleformat{\subsubsection}[runin]{\bfseries}{\bfseries\upshape\sffamily\thesubsubsection.}{1ex}{}

\titlespacing*{\section}{0pt}{*3}{*1}
\titlespacing*{\subsection}{0pt}{*3}{*1}
\titlespacing*{\subsubsection}{0pt}{*2}{*1}

\titlecontents{section}
 [1.5em] %
 {\smallskip}
 {\bfseries\thecontentslabel\hspace{1.02em}}%
{\bfseries}%
 {\,\,\titlerule*[0.77pc]{}\bfseries\contentspage}
\titlecontents{subsection}
 [4em] %
 {\smallskip}
 {\thecontentslabel\hspace{1.02em}}%
{\hspace*{2.32em}}
 {\,\,\titlerule*[0.77pc]{.}\contentspage}

\setlist{topsep=3pt,itemsep=3pt}

\def\wh{\widehat}

\def\lt{\left}
\def\rt{\right}

\def\<{\langle}
\def\>{\rangle}
\def\eps{{\varepsilon}}

\newcommand{\progint}{{\mathrm{int}}}
\newcommand{\progfrac}{{\mathrm{frac}}}
\newcommand{\cH}{{\mathcal{H}}}
\newcommand{\cG}{{\mathcal{G}}}
\newcommand{\cU}{{\mathcal{U}}}
\newcommand{\cW}{{\mathcal{W}}}

\newcommand{\E}{\mathbb{E}}
\newcommand{\cF}{\mathcal{F}}

\author{
Huy Tuan Pham}
\address{\normalfont Department of Mathematics, University of Chicago, Illinois, USA. Research supported by a Clay Research Fellowship, a Sloan Research Fellowship and NSF Award DMS-2543870. Most of this work was completed when the author was supported by a Clay Research Fellowship at the Institute for Advanced Study.}
\email{htpham@uchicago.edu}

\begin{document}
\title{A sharp version of Talagrand's selector process conjecture, with applications to
rounding fractional covers and Bernoulli Sudakov minoration}

\begin{abstract}
    We prove a sharp version of Talagrand's selector process conjecture. Roughly speaking, given any collection of nonnegative weight vectors whose support form a family that is not $p$-small, a random set of density $O(sp)$ captures all but a $2^{-s}$ fraction of weight of some vector with high probability. This gives a common strengthening of Talagrand's selector process conjecture and the Kahn--Kalai conjecture. 
    
    We give two applications of this result. First, towards a conjecture of Talagrand on the equivalence of expectation thresholds and fractional expectation thresholds, we show that a $p$-small fractional cover supported on sets of size at most $t$ can be rounded to a $cp/\log(2t)$-small integral cover. As a corollary, we show that the fractional and integral expectation thresholds are separated by at most a $\log \log$ factor. Second, we prove a Sudakov minoration principle for general positive selector processes, which in particular resolves a problem of Talagrand on Sudakov minoration for the product Bernoulli measure.
\end{abstract}
\maketitle

\section{Introduction}

The study of suprema of stochastic processes is of central interest in probability theory, with influential applications in related areas. While the case if centered Gaussian processes was well-understood thanks to Talagrand's celebrated generic chaining bound and majorizing measure theorem \cite{TalagrandCB,TalagrandCBs,Talagrand1,Talagrand2}, the situation is substantially more delicate beyond the Gaussian setting. In this context, Talagrand has initiated a grand program to understand suprema of general processes, which, thanks to advances on Talagrand's Bernoulli conjecture by Bednorz and Lata\l a \cite{BLatala} and Talagrand's generalized Bernoulli theorem \cite{Talagrand2}, reduces to the understanding of \emph{positive selector processes} as the last missing piece in the picture. In this context, Talagrand's selector process conjecture \cite{Talagrand06, Talagrand, Talagrand2} represents one of the ``Unfulfilled dreams'' \cite{Talagrand2}.  

Among the numerous beautiful conjectures of Talagrand along this program, Talagrand also suggested relations to the study of thresholds of monotone properties and the Kahn--Kalai conjecture on the location of thresholds \cite{KK}. In particular, phenomena around thresholds and positive selector processes were both predicted governed by \emph{first moment obstructions}, which we will introduce precisely later.   

Recent works of Park and the author have confirmed both Talagrand's selector process conjecture \cite{PPT} and the Kahn--Kalai conjecture \cite{PPK}, confirming the first moment obstructions as the fundamental driving force in both cases. While this hints at some connection, and indeed the proof of the Kahn--Kalai conjecture was crucially motivated by the proof of Talagrand's selector process conjecture, no precise connection and implication was earlier known between the two settings. The main result of this paper is a sharp version of Talagrand's selector process conjecture, which gives a common strengthening of Talagrand's selector process conjecture and the Kahn--Kalai conjecture. This gives a precise connection between the study of thresholds of monotone properties in probabilistic combinatorics, and the study of suprema of stochastic processes. As applications of this result, we will show an exponential improvement on the gap between fractional expectation thresholds and expectation thresholds, towards a conjecture of Talagrand \cite{Talagrand}, and show an analog of Sudakov's minoration for positive selector processes, addressing a question of Talagrand.

Given a finite set $X$ and $p\in[0,1]$, let $X_p$ be the random subset of $X$ obtained by retaining each element independently with probability $p$, and denote its law on $2^X$ by $\mu_p$. For $S\subseteq X$, let $\langle S\rangle=\{A\subseteq X:S\subseteq A\}$, and, for $\cG\subseteq2^X$, let $\langle\cG\rangle=\bigcup_{G\in\cG}\langle G\rangle$. We say that $\cG$ is a \emph{cover} of $\cH$ if $\cH \subseteq \langle \cG\rangle$. A family $\cH\subseteq2^X$ is \emph{$p$-small} if there is a cover $\cG\subseteq2^X$ such that 
\[
    \sum_{G\in\cG}p^{|G|}\le\frac12.
\]
We call such $\cG$ a \emph{first moment obstruction} of $\cH$.

Given a collection $\Lambda$ of nonnegative weight vectors on $X$, let
\[
    Z_\lambda(A) = \sum_{x\in A}\lambda(x), \qquad Z_\Lambda(A)=\sup_{\lambda\in\Lambda}Z_\lambda(A).
\]
The random variables $Z_\lambda(X_p)$ form a selector process, and $Z_\Lambda(X_p)$ is its supremum. Talagrand \cite{Talagrand06,Talagrand,Talagrand2} conjectured the following result, which was proved by Park and the author \cite{PPT}.

\begin{theo}[Talagrand's selector process conjecture]\label{thm:tail-formulation}
    There exists an absolute constant $L>0$ such that the following holds. Let $0<p\le1$, and let $\Lambda$ be a collection of nonnegative weight vectors on $X$ such that $\E Z_\Lambda(X_p)<\infty$. Then
    \begin{equation}\label{eq:selector-tail}
        \{A\subseteq X:Z_\Lambda(A)>L\E Z_\Lambda(X_p)\}
        \quad\textrm{is $p$-small}.
    \end{equation}
\end{theo}

Equivalently, the event in (\ref{eq:selector-tail}) is contained in $\langle\cG\rangle$ for a collection $\cG$ satisfying $\sum_{G\in\cG}p^{|G|}\le1/2$. Since $\mu_p(\langle G\rangle)=p^{|G|}$, the sets $\langle G\rangle$ play, for the Bernoulli product measure, the role played by half-spaces in the Gaussian setting.

Positive selector processes represent processes where no cancellation is available, and their suprema form an important part of the non-Gaussian theory. Earlier results \cite{BLatala, Talagrand2} reduce the characterization of suprema of important centered stochastic processes, such as centered Bernoulli-$p$ processes, to positive selector processes. Positive selector processes also arise naturally as discrete approximations of positive empirical processes. %
This connection was used by Park and the author \cite{PPT} to obtain an analogue of the selector process conjecture for general positive empirical processes, answering a question of Talagrand. The original proof of the selector process conjecture was subsequently extended to positive infinitely divisible processes and positive canonical processes by Bednorz, Martynek and Meller \cite{BMM, BMM2}.%

As observed by Talagrand \cite{Talagrand}, the selector process conjecture is equivalent, up to a change in the absolute constant, to the following formulation.
\begin{theo}\label{thm:selector}
    There exists a constant $c>0$ such that the following holds. Let $0<p\le1$, let $\cH\subseteq 2^X$ be not $p$-small, and, for each $H \in \cH$, consider a weight vector $\lambda_{H}:X\to [0,1]$ such that $\lambda_{H}$ is supported on $H$ and $\sum_{x\in H}\lambda_H(x) \ge 1$. Then,
    \[
        \mathbb{E} \max_{H\in \cH} \sum_{x\in X_{p}\cap H} \lambda_H(x) \ge c.
    \]
\end{theo}

Another area where first moment obstructions play a fundamental role is the study of thresholds, which is an important subject in probabilistic combinatorics and random graph theory. Given any nontrivial increasing property $\langle\cH\rangle$ ($\cH$ can be identified as the collection of minimal elements of the property), its threshold $p_c(\cH)$ is the value of $p$ for which
\[
    \P(X_p\in\langle\cH\rangle)=\frac12.
\]

It follows immediately from the union bound that whenever $\cH$ is $p$-small, i.e. a first moment obstruction exists for $\cH$ at density $p$, then $p_c(\cH)\ge p$. Kahn and Kalai \cite{KK} defined the expectation threshold $p_E(\cH)$ as the smallest $p$ such that no first moment obstruction exists, which presents a lower bound for $p_c(\cH)$. 

The Kahn--Kalai conjecture asserts that this elementary lower bound always determines the threshold up to a logarithmic factor. It was resolved by Park and the author \cite{PPK} in the more precise form below.
\begin{theo}\label{thm:kk}
	There exists a constant $C>0$ such that the following holds. For all $\cH$ which is not $p$-small and $\ell = \max_{H\in \cH}|H|$, 
	\[
		\mathbb{P}(X_{Cp\log \ell} \in \langle \cH\rangle) > 1 - o_\ell(1) > 1/2. 
	\]
	In particular, 
	\[
    p_c(\cH)\le C p_E(\cH)\log(2\ell).
\]
\end{theo}

The main result of the present paper is a sharp version of Talagrand's selector process conjecture, Theorem \ref{thm:selector}, which provides a unifying common strengthening of both Theorem \ref{thm:selector} and Theorem \ref{thm:kk}.

\subsection{The sharp selector theorem}

The following theorem is the main result of the paper. 
\begin{theo}\label{thm:sharp-selector}
    There exists a constant $C>0$ such that the following holds. Let $0<p\le1$, let $\cH\subseteq 2^X$ be not $p$-small, and, for each $H \in \cH$, consider a weight vector $\lambda_{H}:X\to [0,1]$ such that $\lambda_{H}$ is supported on $H$ and $\sum_{x\in H}\lambda_H(x) \ge 1$. Let $s$ be a positive integer. Then, with probability at least $1-2^{-s-1}$,
    \[
        \max_{H\in \cH} \sum_{x\in X_{Csp}\cap H} \lambda_H(x) \ge 1-2^{-s}.
    \]
\end{theo}

This result also yields as a directly corollary the following sharp counterpart of the tail-event formulation Theorem \ref{thm:tail-formulation}. 

\begin{corl}\label{cor:sharp-selector-tail}
    There exists a constant $c>0$ such that the following holds. Let $0<p\le1$ and let $\Lambda$ be a collection of nonnegative weight vectors on $X$ such that $\E Z_\Lambda(X_p)<\infty$. For every $0<\delta\le1$, the family
    \[
        \{A\subseteq X:Z_\Lambda(A)>(1+\delta)\E Z_\Lambda(X_p)\}
    \]
    is $q$-small for
    \[
        q=\frac{cp}{\log(2/\delta)}.
    \]
\end{corl}

Thus, in the original tail-event formulation, the far-tail event $Z_\Lambda(A)>L\E Z_\Lambda(X_p)$ can be replaced by the optimal tail event $Z_\Lambda(A)>(1+\delta)\E Z_\Lambda(X_p)$, at the cost of replacing $p$-smallness by $cp/\log(2/\delta)$-smallness. %

Writing $\delta=2^{-s}$, Theorem \ref{thm:sharp-selector} replaces the fixed positive fraction in Theorem \ref{thm:selector} by $1-\delta$, at the cost of increasing the density by a factor of order $\log(1/\delta)$. %

For $s=1$, Theorem \ref{thm:sharp-selector} recovers, up to a change in the sampling density by an absolute constant, the selector process theorem of \cite{PPT}. More importantly, it shows that one can capture all but an $\eps$-fraction of the weight by increasing the density only by a factor of order $\log(1/\eps)$. This dependence is best possible up to an absolute constant, as shown by the example at the beginning of Section \ref{sec:sharp-Tal}. 

Theorem \ref{thm:sharp-selector} also immediately implies the Kahn--Kalai theorem. Indeed, take $\lambda_H(x)=1/|H|$ for $x\in H$ and let $s=\lceil\log(2\ell)\rceil$, where $\ell=\max_{H\in\cH}|H|$. The conclusion gives some $H\in\cH$ with
\[
    |H\setminus X_{Csp}|\le |H|2^{-s}\le\frac12,
\]
and hence $H\subseteq X_{Csp}$. %
Hence, the selector process setup, and the sharp version Theorem \ref{thm:sharp-selector} provides a strengthening of the Kahn--Kalai conjecture in the more general weighted setup. 

While in general the $\log \ell$ factor is needed in the Kahn--Kalai conjecture, Theorem \ref{thm:sharp-selector} describes an essentially optimal unifying picture interpolating the phenomena between $p_E(\cH)$ and $p_c(\cH)\le Cp_E(\cH)\log (2\ell)$.

In the next subsections, we describe several further applications of Theorem \ref{thm:sharp-selector}. 

\subsection{Thresholds and fractional covers}

The expectation threshold can be equivalently expressed via the integer linear program
\begin{align}
    c_{\progint}(\cH; p) := &\min \sum_{W\in 2^X} x_W p^{|W|}\nonumber\\
    &\textrm{subject to }\sum_{W\subseteq H} x_W \ge 1 \qquad \forall H\in \cH,\nonumber\\
    &\qquad \qquad \,\,\,\,\,\,\,\,  x_W\in \{0,1\} \qquad \forall W\in 2^X,\label{eq:int-linprog}
\end{align}

Given the challenge in controlling integer linear programs, Talagrand \cite{Talagrand} introduced the fractional relaxation
\begin{align}
    c_{\progfrac}(\cH; p) := &\min \sum_{W\in 2^X} x_W p^{|W|}\nonumber\\
    &\textrm{subject to }\sum_{W\subseteq H} x_W \ge 1 \qquad \forall H\in \cH,\nonumber\\
    &\qquad \qquad \,\,\,\,\,\,\,\,  x_W\in [0,1] \qquad \,\,\forall W\in 2^X.\label{eq:frac-linprog}
\end{align}
We call $w:2^X\to[0,1]$ satisfying $\sum_{W\subseteq H}w(W)\ge1$ for every $H\in\cH$ a \emph{fractional cover} of $\cH$. The \emph{fractional expectation threshold} $p_f(\cH)$ is the largest $p$ for which $c_{\progfrac}(\cH;p)\le1/2$. We always have
\[
    p_E(\cH)\le p_f(\cH)\le p_c(\cH).
\]
The first inequality follows from the relaxation. For the second, if $w$ is a fractional cover, then
\[
    \mathbb I(X_p\in\langle\cH\rangle)\le\sum_{W\subseteq X_p}w(W),
\]
whose expectation is $\sum_{W\in2^X}w(W)p^{|W|}$.

Talagrand made a bold conjecture that the expectation threshold and fractional expectation threshold are always equivalent up to constant factors. 
\begin{conj}[Talagrand \cite{Talagrand}]\label{conj:Tal-equiv}
    There exists a constant $C>0$ such that
    \[
        p_E(\cH)\le p_f(\cH)\le Cp_E(\cH).
    \]
    Equivalently, given $w:2^X\to[0,1]$ such that $\sum_{W\in2^X}w(W)p^{|W|}\le1/2$, any $\cH$ fractionally covered by $w$ is $p/C$-small.
\end{conj}

Conjecture \ref{conj:Tal-equiv} can be viewed from the perspective of rounding fractional relaxations of integer linear programs. In particular, it provides a broad general class of set-cover-type integer program (\ref{eq:int-linprog}) for which there is no gap between the integral version and its relaxation (\ref{eq:frac-linprog}). %

Our first application of Theorem \ref{thm:sharp-selector} gives the following quantitative form of Talagrand's conjecture for fractional covers with bounded support size.

\begin{theo}\label{thm:main}
    There exists a constant $c>0$ such that the following holds. Let $0<p\le1$ and let $t\ge1$ be an integer. Assume that $\cH$ admits a fractional cover $w : 2^X \to [0,1]$ such that $\sum_{W\in 2^X} w(W)p^{|W|} \le 1/2$ and $w$ is supported on sets of size at most $t$. Then $c_{\progint}(\cH; q) \le 1/2$ for $q = c p / \log(2t)$, i.e. there exists $g: 2^X \to \{0,1\}$ such that $\sum_{W\subseteq H} g(W) \ge 1$ for all $H\in \cH$ and
    \[
        \sum_{W\in2^X} g(W)q^{|W|} \le 1/2.
    \]
\end{theo}

In particular, Theorem \ref{thm:main} verifies the asserted constant-factor rounding for every instance in which $t$ is bounded by an absolute constant. Observe that the factor $\log(2t)$ here is different from the logarithmic factor in the Kahn--Kalai theorem. Even for bounded $t$, the sets in $\cH$ may have arbitrarily large size.

The proof of Theorem \ref{thm:main} provides a novel mechanism for rounding fractional covers that is significantly different from standard randomized rounding, where the key driving force is Theorem \ref{thm:sharp-selector}. %

Combining Theorem \ref{thm:main} with the sampling argument of Fischer and Person \cite{FP} gives the universal estimate
\[
    p_f(\cH)\le C p_E(\cH)\log\log(4|X|).
\]
Prior to our work, the best known bound on the gap between $p_f(\cH)$ and $p_E(\cH)$ is $p_E(\cH)\le p_f(\cH)\le O(p_E(\cH)\log |X|$, which follows directly from the Kahn--Kalai conjecture. 
\begin{corl}\label{corl:loglog-gap}
    There exists an absolute constant $c>0$ such that the following holds. Let $X$ be a nonempty finite set, let $0<p\le1$, and suppose that $\cH\subseteq2^X$ admits a fractional cover $w:2^X\to[0,1]$ such that
    \[
        \sum_{W\in2^X}w(W)p^{|W|}\le\frac12.
    \]
    Then $\cH$ is $q$-small for
    \[
        q=\frac{cp}{\log\log(4|X|)}.
    \]
    Consequently, for every nontrivial increasing family $\langle\cH\rangle\subseteq2^X$,
    \[
        p_E(\cH)\le p_f(\cH)\le C p_E(\cH)\log\log(4|X|)
    \]
    for an absolute constant $C>0$.
\end{corl}
In particular, Corollary \ref{corl:loglog-gap} provides an exponential improvement over the previous bound, and strictly separates the fractional expectation threshold $p_f(\cH)$ from the threshold $p_c(\cH)$. 

Previous progress on Conjecture \ref{conj:Tal-equiv} was restricted to several special cases. Talagrand \cite{Talagrand} proved the case where $w$ is supported on singletons. Frankston, Kahn and Park \cite{FKP} proved the case where $w$ is supported on pairs. Fischer and Person \cite{FP} considered fractional covers supported on nearly linear uniform hypergraphs, and DeMarco and Kahn \cite{DK} treated the case of cliques in graphs. The present work is independent of the work of Dubroff, Kahn and Park \cite{DKP}, who, using a quantitative form of the selector theorem in \cite{BMM}, obtained a weaker version of Theorem \ref{thm:main} with polynomial dependence on $t$, and hence does not provide improvement over the $\log |X|$ bound for the gap between $p_f(\cH)$ and $p_E(\cH)$ in Corollary \ref{corl:loglog-gap}. %

\subsection{Sudakov minoration for $\mu_p$}

Sudakov minoration is one of the basic lower-bound principles for Gaussian processes. In one of its standard forms, if $(G_t)_{t\in T}$ is a centered Gaussian process and $T$ contains $N$ points whose pairwise canonical distances are at least $a$, then
\[
    \E\sup_{t\in T}G_t\ge c a\sqrt{\log N}
\]
for an absolute constant $c>0$. This is the fundamental one-scale entropy lower bound and an important part of the majorizing measure theory; see Sudakov's original work \cite{Sudakov71} and the accounts in \cite{LedouxTalagrand,Talagrand1,Talagrand2}.

Given its central importance in the theory of Gaussian processes, Talagrand asked for an analogue of Sudakov minoration for the Bernoulli product measure $\mu_p$. In \cite[Problem 4.1]{TalFav}, he proposed the following concrete formulation.

\begin{prob}[Talagrand's Problem 4.1]\label{prob:talagrand-sudakov}
    Does there exist an absolute constant $c>0$ such that the following holds? Let $0<p<1$, let $X$ be a finite set, let $n\ge1$, and let $\cF\subseteq\binom{X}{n}$ satisfy
    \[
        |\cF|\ge p^{-n}
        \qquad\textrm{and}\qquad
        |F\cap F'|\le n/2
        \quad\textrm{for all distinct }F,F'\in\cF.
    \]
    Must one have
    \[
        \E\max_{F\in\cF}|F\cap X_p|\ge cn?
    \]
\end{prob}

The overlap assumption plays the role of separation.%

We will prove the following more general weighted form of Sudakov minoration for selector processes that immediately implies an answer to Problem \ref{prob:talagrand-sudakov}. For a nonnegative weight vector $\lambda$ and $W\subseteq X$, we define $\lambda(W)=\sum_{x\in W}\lambda(x)$.

\begin{theo}[Sudakov minoration for selector processes]\label{thm:sudakov}
    There exists an absolute constant $C>0$ such that the following holds. Let $X$ be finite, let $0<p<1$, let $n$ and $s$ be positive integers, let $r>0$, and let $\Lambda_0$ be a finite subcollection of $\Lambda\subseteq\R_{\ge0}^X$. Suppose that, for every $\lambda\in\Lambda_0$, there is a set $C_\lambda\subseteq X$ with $|C_\lambda|\le n$, and that
    \[
        |\Lambda_0|\ge p^{-n}.
    \]
    If, for all distinct $\lambda,\lambda'\in\Lambda_0$,
    \begin{equation}\label{eq:sudakov-private-mass}
        \lambda(C_\lambda\setminus C_{\lambda'})\ge r,
    \end{equation}
    then, with probability at least $1-2^{-s-1}$,
    \[
        Z_\Lambda(X_{Csp})\ge(1-2^{-s})r.
    \]
\end{theo}

In particular, under the remaining hypotheses of Theorem \ref{thm:sudakov}, for $0<\eps\le1/2$ we have
\begin{equation}\label{eq:sudakov-fixed-density}
    |\Lambda_0|\ge\left(\frac{C\log(2/\eps)}{p}\right)^n
    \quad\Longrightarrow\quad
    \P\bigl(Z_\Lambda(X_p)\ge(1-\eps)r\bigr)\ge1-\frac{\eps}{2}. 
\end{equation}

The expected value consequence also follows immediately. 

\begin{corl}\label{corl:sudakov-expectation}
    There exists an absolute constant $c>0$ such that the following holds. Let $X$ be finite, let $0<p<1$, let $n$ be a positive integer, let $r>0$, and let $\Lambda_0$ be a finite subcollection of $\Lambda\subseteq\R_{\ge0}^X$. Suppose that there are sets $C_\lambda\subseteq X$, $|C_\lambda|\le n$, such that $|\Lambda_0|\ge p^{-n}$ and (\ref{eq:sudakov-private-mass}) holds for all distinct $\lambda,\lambda'\in\Lambda_0$. Then
    \[
        \E Z_\Lambda(X_p)\ge cr.
    \]
\end{corl}

Specializing to indicator vectors gives a positive answer to Problem \ref{prob:talagrand-sudakov}.

\begin{corl}[Talagrand's $\mu_p$-Sudakov minoration]\label{corl:talagrand-sudakov}
    There exists an absolute constant $c>0$ such that the following holds. Let $0<p\le1$, let $X$ be a finite set, let $n\ge1$, and let $\cF\subseteq\binom{X}{n}$ satisfy
    \[
        |\cF|\ge p^{-n}
        \qquad\textrm{and}\qquad
        |F\cap F'|\le n/2
        \quad\textrm{for all distinct }F,F'\in\cF.
    \]
    Then
    \[
        \E\max_{F\in\cF}|F\cap X_p|\ge cn.
    \]
\end{corl}

Theorem \ref{thm:sudakov} is a direct discrete analogue of Sudakov minoration. %
The proof of this result uses crucially Theorem \ref{thm:sharp-selector}. %

\subsection{Organization}

In Section \ref{sec:sharp-Tal}, we prove the main result, Theorem \ref{thm:sharp-selector}, via the new notion of \emph{towers of minimum fragments}. 

In Section \ref{sec:fractional-covers}, we prove Theorem \ref{thm:main} and Corollary \ref{corl:loglog-gap}. In Section \ref{sec:sudakov}, we prove Theorem \ref{thm:sudakov} and its consequences using Theorem \ref{thm:sharp-selector}.

We make no attempt to optimize absolute constants. Throughout the paper, $\mathbb N=\{0,1,2,\ldots\}$ and logarithms are in base $2$. Whenever a displayed density is at least $1$, we interpret the corresponding sample as the full set $X$. 

A preliminary version of the selector theorem and its application to fractional covers appeared in the extended abstract \cite{PhamSharp}. The present journal version follows similar key ideas in defining the \emph{towers of minimum fragments}, while streamlining the proof and strengthening the failure probability in the main result.

\section{A sharp version of Talagrand's selector process conjecture}\label{sec:sharp-Tal}

In this section, we prove the sharp version of Talagrand's selector process conjecture, Theorem \ref{thm:sharp-selector}. 

Before moving on to the proof, we first remark that Theorem \ref{thm:sharp-selector} is sharp in that, to guarantee that $\max_{H\in \cH} \sum_{x\in H\cap X_q} \lambda_H(x) \ge 1 - \eps$, we need in general that $q \ge \Omega(p \log (1/\eps))$. Indeed, consider the example of the Erd\H{o}s-R\'enyi random graph $G(n,q)$ (so $X = \binom{[n]}{2}$), where $n$ is even, and let $\cH$ be the collection of perfect matchings, which is not $c/n$-small for an absolute constant $c>0$. For $q = s/n$ and $ne^{-s}$ sufficiently large, $G(n,q)$ contains at least $\Omega(e^{-s}n)$ isolated vertices with high probability. Hence, for $\lambda_H(x) = 1/|H|$ for each $x\in H$ and $H\in \cH$, we have that $\max_{H\in \cH} \sum_{x\in H\cap X_q} \lambda_H(x) \le 1 - \Omega(e^{-s})$.  

\subsection{Proof strategy and difficulty}
To motivate the overall strategy of the proof, we will first give a high level discussion of the plan in \cite{PPT} and explain why proving Theorem \ref{thm:sharp-selector} requires a new idea. 

For $S\subseteq X$ and $\lambda_H:X\to [0,1]$ supported on $H$, we denote $\lambda_{H}(S) = \sum_{x\in S\cap H}\lambda_H(x)$. The proof of Talagrand's selector process conjecture in \cite{PPT} relies on the key notion of \emph{minimum fragments}. Roughly speaking, for every $W\in 2^X$ and $H \in \cH$, the minimum fragment $T(W,H)$ is a subset of $H$ that arises from its interaction with $W$. One then observes that the collection $\cU(W) := \{T(W,H) : H\in \cH\}$ forms a cover for $\cH$ for each $W\in 2^X$. Furthermore, under the correct definition of the minimum fragment, one can show a key property that allows the fragments to be efficiently encoded. Defining $W$ to be bad if $\max_{H\in \cH} \lambda_H(W) < c$, one then uses this property to prove, for $q=Cp$,
\begin{equation}\label{eq:bad-lowcost-old}
    \E_{W\sim X_q}\lt[\mathbb{I}(W \textrm{ is bad}) \sum_{U\in \cU(W)} p^{|U|}\rt] \ll 1/2,
\end{equation}
which leads to an upper bound on the probability that $X_q$ is bad, since any cover $\cU$ for $\cH$ must have $\sum_{U\in \cU}p^{|U|} \ge 1/2$. 

Our proof of Theorem \ref{thm:sharp-selector} follows the same overarching scheme. However, the notion of minimum fragments as defined in \cite{PPT} cannot be used directly, and we need a more involved notion, \emph{towers of minimum fragments}. Without recalling the detailed definition of minimum fragments in \cite{PPT}, we briefly explain why they are not sufficient to prove Theorem \ref{thm:sharp-selector}. 

The selector process conjecture yields, for a typical sample $W_1 \sim X_{q}$, a heavy set $\tilde{H} \in \cH$ for which $\lambda_{\tilde{H}}(W_1) \ge c$. One may hope to obtain Theorem \ref{thm:sharp-selector} by replacing $\cH$ by the collection of sets $\tilde{H}\setminus W_1$ and iterating with new samples $W_i\sim X_q$. This is roughly the iteration behind the proof of the Kahn--Kalai conjecture. However, in the general weighted setup of the selector process conjecture, there is a fundamental obstruction to this iteration. Roughly speaking, based on each fragment $T(W_1,H)$, one obtains a heavy set $\tilde{H}\in \cH$ depending appropriately on $W_1$ and $H$. The proof requires $\tilde{H}$ to be allowed to use elements outside $W_1 \cup H$. Hence, we lose any control on the $p$-smallness of the collection of sets of the form $\tilde{H} \setminus W_1$. 

Instead, our proof of Theorem \ref{thm:sharp-selector} defines the notion of towers of minimum fragments. For each $H\in \cH$ and samples $W_1,\dots,W_\ell\sim X_q$, where $\ell=O(s)$, the tower of minimum fragments is a suitable collection of disjoint subsets $T_1,\dots,T_\ell$ of $H$. It combines the information from the interaction of $H$ with all the samples $W_1,\dots,W_\ell$ across the $\ell$ ``iterations''. In the context of Theorem \ref{thm:sharp-selector}, say that $(W_1,\dots,W_\ell)$ is bad if $\max_{H\in \cH} \lambda_H(\bigcup_{i=1}^{\ell}W_i) < 1-2^{-s}$. The towers of minimum fragments are defined so that they satisfy a key property allowing them to be efficiently encoded. As a result, for any bad $(W_1,\dots,W_\ell)$, we obtain an analogue of (\ref{eq:bad-lowcost-old}). 

Achieving the correct notion of the tower of minimum fragments is a significant challenge. The design of the tower of minimum fragments is guided by the following principle:
\begin{itemize}
    \item $(T_1,\dots,T_\ell)$ is defined so that $(\mathbf Z,\mathbf t)$, where $Z_i=W_i \cup T_i$ and $t_i=|T_i|$, satisfies a suitable feasibility condition with a witness $(r,H')$, where $H'\in\cH$. 
    \item For a feasible tuple $(\mathbf Z,\mathbf t)$, we can guarantee that $\max_{H \in \cH} \lambda_H(\bigcup_{i=1}^{\ell} Z_i) \ge 1-2^{-s}$. 
    \item The fragments $T_i$ can be identified efficiently in a subset $\wh T_i$ which is uniquely determined by the feasible data $(\mathbf Z,\mathbf t)$. 
    \item The size $|\wh T_i|$ is comparable to $|T_i|$. 
\end{itemize}
Achieving all properties above simultaneously is significantly harder than the case $s=1$ considered previously in \cite{PPT, BMM, BMM2}. In particular, a key idea there is to verify that if $T_i$ has elements outside the desired set $\wh T_i$, then a smaller fragment can be constructed. In our case, this is highly nontrivial since feasibility has to be defined jointly across $\ell$ iterations. Shrinking a fragment $T_i$ at an earlier step leaves additional elements uncovered in the later iterations, which might not be covered by any modification of the later fragments. Thus, our construction of the tower of minimum fragments is significantly different from the previous constructions.

In the next three subsections, we first make some simple reductions and then formally introduce towers of minimum fragments together with their key property. In Subsection \ref{subsec:proof}, we complete the proof of Theorem \ref{thm:sharp-selector}.

\subsection{Simple reductions}

We first reduce Theorem \ref{thm:sharp-selector} to a statement for weight vectors taking positive values among $\{B^{-z}:z\in\mathbb{N}\}$. 

\begin{theo}\label{thm:coarse-selector}
    There exists a constant $C_0>0$ such that the following holds. Let $0<p\le1$, let $B=2^L$ for an integer $L\ge 2$, and let $\cH\subseteq 2^X$ be not $p$-small. For each $H \in \cH$, consider a weight vector $\lambda_{H}:X\to [0,1]$ such that $\lambda_{H}$ is supported on $H$, 
    \[
        \frac{1}{B}\le \lambda_H(H)\le 1,
    \]
    and every positive value of $\lambda_H$ belongs to $\{B^{-z}:z\in\mathbb{N}\}$. Then, with probability at least $1-2B^{-3}$, there exists $H\in\cH$ such that
    \[
        \lambda_H(H\setminus X_{C_0pL})\le B^{-2}.
    \]
\end{theo}

We show that Theorem \ref{thm:coarse-selector} implies Theorem \ref{thm:sharp-selector}. 
\begin{proof}
    For each $H\in \cH$, let $\nu_H=\lambda_H/\lambda_H(H)$, so $\nu_H(H)=1$. For $\nu_H(x)>0$, let $\lambda'_H(x)=B^{\lfloor \log_B\nu_H(x)\rfloor}$, and let $\lambda'_H(x)=0$ when $\nu_H(x)=0$. Then
    \[
        \lambda'_H(x)\le \nu_H(x)\le B\lambda'_H(x).
    \]
    In particular, $B^{-1}\le\lambda'_H(H)\le 1$, and every positive value of $\lambda'_H$ belongs to $\{B^{-z}:z\in\mathbb{N}\}$. Furthermore, for every $Y\subseteq X$, 
    \[
        \lambda'_H(H\setminus Y)\le B^{-2}\qquad\Longrightarrow\qquad \lambda_H(Y)\ge1-B^{-1}.
    \]
    Indeed, in this case,
    \[
        \nu_H(H\setminus Y)\le B\lambda'_H(H\setminus Y)\le B^{-1},
    \]
    and $\lambda_H(Y)=\lambda_H(H)\nu_H(Y)\ge \nu_H(Y)\ge1-B^{-1}$.

    We now take $B=4^s=2^{2s}$ and apply Theorem \ref{thm:coarse-selector} with $L=2s$ to the weight vectors $\lambda'_H$. Since $B^{-1}=4^{-s}\le 2^{-s}$, Theorem \ref{thm:sharp-selector} follows with $C=2C_0$. Moreover, its failure probability is at most
    \[
        2B^{-3}=2^{1-6s}\le2^{-s-1}.
    \]
\end{proof}

\subsection{Towers of minimum fragments}

Fix $B=2^L$ for an integer $L\ge 2$, and let $\ell=3L+2$.
Notice that
\begin{equation}\label{eq:ell-choice}
    2^{-\ell}B=\frac{1}{4}B^{-2}.
\end{equation}

Replacing $\cH$ by $\{\operatorname{supp}(\lambda_H):H\in\cH\}$ and retaining one associated weight vector for each resulting set preserves the fact that the family is not $p$-small and does not change the conclusion. As such, without loss of generality, we may assume
\[
    H=\operatorname{supp}(\lambda_H),
\]
for every $H\in\cH$. In particular, $\lambda_H(x)>0$ for every $x\in H$.

For $H\in \cH$ and $j\in\mathbb{N}$, let
\[
    H_j:=\{x\in H:\lambda_H(x)=B^{-j}\}.
\]
Let $M_H=\max\{j:H_j\ne\emptyset\}$, and put $H_j=\emptyset$ for all other integers $j$. We denote $H_{<r}=\bigcup_{j<r}H_j$ and $H_{\ge r}=\bigcup_{j\ge r}H_j$. For $Y\subseteq X$, we also denote
\[
    \lambda_{\ge r}(Y)=\sum_{x\in Y\cap H_{\ge r}}\lambda_H(x).
\]
All level sets and weight notation associated with a displayed $H$ are formed using $\lambda_H$. In particular, we suppress $\lambda_H$ from the arguments of $\Phi$ and its variants.

For $r\in \{0,\dots,M_H+3\}$, define
\begin{equation}\label{eq:def-Phi}
    \Phi(r;Y,H):=\lambda_{\ge r}(Y)-(1-B^{-2})\lambda_{\ge r}(H)+B^{-r}\left(|H_{<r-2}|+\frac{3}{4}|H_{r-2}|+\frac{3}{4B}|H_{r-1}|\right).
\end{equation}

For an integer $n\ge 0$, define
\[
    \lceil n\rceil_B=\begin{cases}
        0,&n=0,\\
        B^{\lceil\log_Bn\rceil},&n\ge 1.
    \end{cases}
\]
Thus, for $n\ge 1$,
\begin{equation}\label{eq:B-adic-envelope}
    n\le \lceil n\rceil_B<Bn.
\end{equation}

Given sets $\mathbf{Z}=(Z_1,\dots,Z_\ell)$ and nonnegative integers $\mathbf{t}=(t_1,\dots,t_\ell)$, let
\[
    \overline{Z}_i=\bigcup_{h\le i}Z_h\qquad (i\in[\ell])
\]
and $\overline{Z}_0=\emptyset$. Define
\begin{equation}\label{eq:constraint-1}
    \overline{\Phi}(r;\mathbf{Z},\mathbf{t},H):=\Phi(r;\overline{Z}_\ell,H)-B^{-r}\sum_{h=1}^{\ell}t_h.
\end{equation}

\begin{defn}\label{def:feas}
    We say that $(\mathbf{Z},\mathbf{t})$ is \emph{feasible} if there exists $H\in \cH$ and an integer $r\in\{0,\dots,M_H+3\}$ such that
    \begin{equation}\label{eq:constraint-2}
        H_{<r-2}\subseteq Z_1,
    \end{equation}
    and, for every $j\in\{1,2\}$ and $i\in[\ell]$, at least one of the following holds:
    \begin{align}
        |H_{r-j}\setminus\overline{Z}_i|+\sum_{h\le i}t_h&\le 2^{-i}\lceil |H_{r-j}|\rceil_B;\label{eq:constraint-3}\\
        H_{r-j}\setminus\overline{Z}_{i-1}&\subseteq Z_i.\label{eq:constraint-4}
    \end{align}
    Furthermore,
    \begin{equation}\label{eq:constraint-5}
        \overline{\Phi}(r;\mathbf{Z},\mathbf{t},H)\ge 0.
    \end{equation}
    We refer to $(r,H)$ as a witness of the feasibility of $(\mathbf{Z},\mathbf{t})$.
\end{defn}

Define the lexicographic ordering on $\mathbb{Z}^{\ell}$ by $\mathbf{x}\preceq\mathbf{y}$ if $\mathbf{x}=\mathbf{y}$ or there exists $i\in[\ell]$ such that $x_i<y_i$ and $x_j=y_j$ for all $j<i$, where $\mathbf{x}=(x_1,\dots,x_\ell)$ and $\mathbf{y}=(y_1,\dots,y_\ell)$.
\begin{defn}[Towers of minimum fragments]
    Given a tuple $\mathbf{W}=(W_1,\dots,W_\ell)$ of subsets of $X$ and $H\in\cH$, we say that a tuple of sets $\mathbf{T}=(T_1,\dots,T_\ell)$ is a \emph{tower of fragments} of $(\mathbf{W},H)$ if the following holds.
    \begin{itemize}
        \item The sets $T_1,\dots,T_\ell$ are pairwise disjoint subsets of $H$, and, for every $i\in[\ell]$,
        \[
            T_i\cap\left(\bigcup_{h\le i}W_h\cup\bigcup_{h<i}T_h\right)=\emptyset.
        \]
        \item The tuple $\mathbf{Z}=\mathbf{Z}(\mathbf{W},\mathbf{T})=(Z_1,\dots,Z_\ell)$, where $Z_i=W_i\cup T_i$, together with $\mathbf{t}=(|T_1|,\dots,|T_\ell|)$ is feasible.
    \end{itemize}

    We say that $\mathbf{T}=\mathbf{T}(\mathbf{W},H)$ is the \emph{minimum tower of fragments} of $(\mathbf{W},H)$ if $(|T_1|,\dots,|T_\ell|)$ is minimal among all towers of fragments of $(\mathbf{W},H)$ in the lexicographic ordering. We fix an arbitrary deterministic choice among towers with the same size vector.
\end{defn}

\subsection{Key properties of the tower of minimum fragments}
In this subsection, we record several key properties of the tower of minimum fragments.

Let $\mathbf{T}=\mathbf{T}(\mathbf{W},H)$ be a minimum tower of fragments of $(\mathbf{W},H)$, and let $\mathbf{Z}=\mathbf{Z}(\mathbf{W},\mathbf{T})$ be defined by $Z_i=W_i\cup T_i$. Let $\mathbf{t}=(|T_1|,\dots,|T_\ell|)$.

We first record some simple properties that follow directly from the definition. 

\begin{lemma}
    For every $(\mathbf{W},H)$, a tower of fragments exists.     
\end{lemma}
\begin{proof}
    Consider $T_1=H\setminus W_1$ and $T_i=\emptyset$ for $i>1$. We use the witness $H$ and the threshold $r=M_H+3$. In this case, $H_{\ge r}=H_{r-1}=H_{r-2}=\emptyset$ and $H_{<r-2}=H\subseteq Z_1$. Furthermore,
    \[
        \overline{\Phi}(r;\mathbf{Z},\mathbf{t},H)=B^{-r}|H|-B^{-r}|T_1|=B^{-r}|H\cap W_1|\ge 0.
    \]
    Thus, all the conditions in Definition \ref{def:feas} hold.
\end{proof}

In particular, this implies the existence of a tower of minimum fragments. 

We first record the consequence of the feasibility conditions for the two levels $H_{r-1}$ and $H_{r-2}$.
\begin{lemma}\label{lem:terminal-remainder}
    If $(\mathbf{Z},\mathbf{t})$ is feasible with witness $(r,H)$, then, for $j\in\{1,2\}$,
    \[
        |H_{r-j}\setminus\overline{Z}_\ell|\le \frac{1}{4}B^{-2}|H_{r-j}|.
    \]
\end{lemma}
\begin{proof}
    At the last step, either (\ref{eq:constraint-3}) holds, in which case
    \[
        |H_{r-j}\setminus\overline{Z}_\ell|+\sum_{h=1}^{\ell}t_h\le 2^{-\ell}\lceil |H_{r-j}|\rceil_B,
    \]
    or (\ref{eq:constraint-4}) holds, in which case $H_{r-j}\setminus\overline{Z}_{\ell-1}\subseteq Z_\ell$ and $H_{r-j}\setminus\overline{Z}_\ell=\emptyset$. If $|H_{r-j}|>0$, (\ref{eq:B-adic-envelope}) and (\ref{eq:ell-choice}) give
    \[
        |H_{r-j}\setminus\overline{Z}_\ell|\le 2^{-\ell}\lceil |H_{r-j}|\rceil_B<2^{-\ell}B|H_{r-j}|=\frac{1}{4}B^{-2}|H_{r-j}|.
    \]
    The statement is immediate when $|H_{r-j}|=0$.
\end{proof}

The next lemma explains the choice of the potential $\Phi$.
\begin{lemma}\label{lem:feasible-capture}
    If $(\mathbf{Z},\mathbf{t})$ is feasible with witness $(r,H)$, then
    \[
        \lambda_H(H\setminus\overline{Z}_\ell)\le B^{-2}\lambda_H(H)\le B^{-2}.
    \]
\end{lemma}
\begin{proof}
    By Lemma \ref{lem:terminal-remainder}, for $j\in\{1,2\}$,
    \begin{align*}
        \lambda_H(\overline{Z}_\ell\cap H_{r-j})-(1-B^{-2})\lambda_H(H_{r-j})
        &=B^{-2}B^{-(r-j)}|H_{r-j}|-B^{-(r-j)}|H_{r-j}\setminus\overline{Z}_\ell|\\
        &\ge \frac{3}{4}B^{-2}\lambda_H(H_{r-j}).
    \end{align*}
    Every element of $H_{<r-2}$ belongs to $Z_1$ and has weight at least $B^{3-r}$. Hence,
    \[
        B^{-2}\lambda_H(H_{<r-2})\ge B^{1-r}|H_{<r-2}|\ge B^{-r}|H_{<r-2}|.
    \]
    Recalling (\ref{eq:def-Phi}), (\ref{eq:constraint-1}), and (\ref{eq:constraint-5}), we obtain
    \begin{align*}
        \lambda_H(\overline{Z}_\ell\cap H)-(1-B^{-2})\lambda_H(H)
        &\ge \Phi(r;\overline{Z}_\ell,H)\\
        &=\overline{\Phi}(r;\mathbf{Z},\mathbf{t},H)+B^{-r}\sum_{h=1}^{\ell}t_h\ge 0.
    \end{align*}
    The left-hand side is $B^{-2}\lambda_H(H)-\lambda_H(H\setminus\overline{Z}_\ell)$, which proves the lemma.
\end{proof}

The key property of towers of minimum fragments roughly shows that each $T_i$ is contained in a set $\wh T_i$ which is determined by $(\mathbf{Z},\mathbf{t})$, while the total size of the sets $\wh T_i$ is comparable to the total size of the fragments. We next fix a choice of a witness of the feasibility of $(\mathbf{Z},\mathbf{t})$ which only depends on $(\mathbf{Z},\mathbf{t})$.

Fix an arbitrary ordering on $\cH$. For each $H\in\cH$, consider the thresholds $r\in\{0,\dots,M_H+3\}$ satisfying (\ref{eq:constraint-2}) and at least one of (\ref{eq:constraint-3}) and (\ref{eq:constraint-4}) for every $j\in\{1,2\}$ and $i\in[\ell]$. If there is no such threshold, we ignore $H$. Otherwise, choose one which maximizes $\overline{\Phi}(r;\mathbf{Z},\mathbf{t},H)$, breaking ties by taking the smallest $r$. Since $(\mathbf{Z},\mathbf{t})$ is feasible, the resulting maximum is nonnegative for at least one $H$. Let $\wh H$ be the first such $H$ and let $\wh r$ be the corresponding threshold. At every pair $(j,i)$, we use (\ref{eq:constraint-3}) whenever it holds and otherwise use (\ref{eq:constraint-4}).

We define
\begin{align*}
    \wh T_1&=\wh H_{<\wh r-2}\cup\bigcup_{\substack{j\in\{1,2\}:\\(\ref{eq:constraint-3})\text{ does not hold for }(j,1)}}\wh H_{\wh r-j},\\
    \wh T_i&=\bigcup_{\substack{j\in\{1,2\}:\\(\ref{eq:constraint-3})\text{ does not hold for }(j,i)}}\left(\wh H_{\wh r-j}\setminus\overline Z_{i-1}\right)\qquad (2\le i\le\ell).
\end{align*}
Here (\ref{eq:constraint-3}) is evaluated with the pair $(\wh r,\wh H)$. 
The sets $\wh T_i$ depend only on $(\mathbf{Z},\mathbf{t})$. Furthermore, for each $j\in\{1,2\}$, there is at most one $i$ for which (\ref{eq:constraint-3}) does not hold and $\wh H_{\wh r-j}\setminus\overline Z_{i-1}$ is nonempty.

The key property of towers of minimum fragments says that $\wh T_i$ always contains the fragment $T_i$. 

\begin{lemma}[Key property of towers of minimum fragments]\label{lem:fragment-contain}
    We have $T_i \subseteq \wh T_i$ for every $i\in[\ell]$. 
\end{lemma}

\begin{proof}[Proof of Lemma \ref{lem:fragment-contain}]
    Suppose that $x\in T_i\setminus\wh T_i$. We first have $x\notin\wh H_{<\wh r-2}$. Indeed, $\wh H_{<\wh r-2}\subseteq Z_1$, which contradicts the definition of a tower of fragments if $i>1$, while $\wh H_{<\wh r-2}\subseteq\wh T_1$ if $i=1$.

    Define
    \[
        T'_h=\begin{cases}
            T_i\setminus\{x\},&h=i,\\
            T_{i+1}\cup\{x\},&h=i+1,\ i<\ell,\ x\notin W_{i+1},\\
            T_h,&\text{otherwise},
        \end{cases}
    \]
    and let $\mathbf T'=(T'_1,\dots,T'_\ell)$, $t'_h=|T'_h|$, $Z'_h=W_h\cup T'_h$, $\mathbf t'=(t'_1,\dots,t'_\ell)$, and $\mathbf Z'=(Z'_1,\dots,Z'_\ell)$. The sets $T'_1,\dots,T'_\ell$ satisfy the first condition in the definition of a tower of fragments, and $\mathbf t' \prec \mathbf t$.

    We show that $(\mathbf Z',\mathbf t')$ is feasible with witness $(\wh r,\wh H)$. Condition (\ref{eq:constraint-2}) follows from $x\notin\wh H_{<\wh r-2}$. Suppose that $x\in\wh H_{\wh r-j}$ for some $j\in\{1,2\}$. Since $x\in T_i$, we have $x\notin\overline Z_{i-1}$. Since $x\notin\wh T_i$, condition (\ref{eq:constraint-3}) holds at $(j,i)$, and
    \[
        |\wh H_{\wh r-j}\setminus\overline Z'_i|+\sum_{h\le i}t'_h
        =|\wh H_{\wh r-j}\setminus\overline Z_i|+\sum_{h\le i}t_h.
    \]
    If $i<\ell$, then $x\in Z'_{i+1}$ and $\overline Z'_h=\overline Z_h$ for every $h\ge i+1$. Thus, at step $i+1$, condition (\ref{eq:constraint-4}) remains valid if it held before, while the left-hand side of (\ref{eq:constraint-3}) does not increase. The same is true at every later step. For the other value of $j$, and also when $x\notin\wh H_{\wh r-1}\cup\wh H_{\wh r-2}$, the relevant set differences are unchanged, while the partial sums in (\ref{eq:constraint-3}) do not increase. Thus, the previously valid condition remains valid. 

    We next consider (\ref{eq:constraint-5}). If $i<\ell$ and $x\notin W_{i+1}$, then $\overline Z'_\ell=\overline Z_\ell$ and $\sum_{h=1}^{\ell}t'_h=\sum_{h=1}^{\ell}t_h$, so $\overline\Phi$ is unchanged. If $i<\ell$ and $x\in W_{i+1}$, then $\overline Z'_\ell=\overline Z_\ell$ and $\sum_{h=1}^{\ell}t'_h=\sum_{h=1}^{\ell}t_h-1$, so $\overline\Phi$ increases by $B^{-\wh r}$. Finally, suppose that $i=\ell$. Then $\overline Z'_\ell=\overline Z_\ell\setminus\{x\}$ and $\sum_{h=1}^{\ell}t'_h=\sum_{h=1}^{\ell}t_h-1$. In this case, $\overline\Phi$ increases by $B^{-\wh r}$ if $x\in\wh H_{\wh r-1}\cup\wh H_{\wh r-2}$ or $x\notin\wh H$, and changes by
    \[
        B^{-\wh r}-\lambda_{\wh H}(x)\ge0
    \]
    if $x\in\wh H_{\ge\wh r}$. Thus, (\ref{eq:constraint-5}) holds.
    
    Hence, $\mathbf T'$ is a tower of fragments with a lexicographically smaller size vector, contradicting the minimality of $\mathbf T$.
\end{proof}

In order for the key property of towers of minimum fragments to be useful, we need the following bound on the sets $\wh T_i$.

\begin{lemma}\label{lem:fragment-contain-2}
    We have
    \[
        \sum_{i=1}^{\ell}|\wh T_i|<6\sum_{i=1}^{\ell}t_i
    \]
    whenever $\sum_{i=1}^{\ell}t_i>0$, and $\wh T_i=\emptyset$ for all $i$ when $\sum_{i=1}^{\ell}t_i=0$.
\end{lemma}
\begin{proof}
    We first claim that
    \begin{equation}\label{eq:large-levels}
        |\wh H_{<\wh r-2}|\le \frac{3}{2}\sum_{i=1}^{\ell}t_i.
    \end{equation}
    The claim is immediate if $\wh H_{<\wh r-2}=\emptyset$. Suppose then that $\wh H_{<\wh r-2}\ne\emptyset$. In particular, the threshold $\wh r-1$ belongs to the allowed range.

    The pair $(\wh r-1,\wh H)$ satisfies (\ref{eq:constraint-2})--(\ref{eq:constraint-4}). Indeed, the conditions on $\wh H_{\wh r-2}$ are the same as before. Also, $\wh H_{\wh r-3}\subseteq\wh H_{<\wh r-2}\subseteq Z_1$, so (\ref{eq:constraint-4}) holds at the first step and the relevant set difference is empty afterwards. Finally, $\wh H_{<\wh r-3}\subseteq Z_1$.

    A direct calculation from (\ref{eq:def-Phi}) gives
    \begin{align}
        &\overline{\Phi}(\wh r-1;\mathbf Z,\mathbf t,\wh H)-\overline{\Phi}(\wh r;\mathbf Z,\mathbf t,\wh H)\nonumber\\
        &\qquad=B^{-\wh r-1}\left(\frac{1}{4}|\wh H_{\wh r-1}|-B^2|\wh H_{\wh r-1}\setminus\overline Z_\ell|\right)\nonumber\\
        &\qquad\qquad+B^{-\wh r}\left(\frac{3B}{4}-1\right)|\wh H_{\wh r-3}|+B^{-\wh r}(B-1)|\wh H_{<\wh r-3}|\nonumber\\
        &\qquad\qquad-B^{-\wh r}(B-1)\sum_{i=1}^{\ell}t_i.\label{eq:shift-Phi}
    \end{align}
    By Lemma \ref{lem:terminal-remainder},
    \[
        |\wh H_{\wh r-1}\setminus\overline Z_\ell|\le\frac{1}{4}B^{-2}|\wh H_{\wh r-1}|,
    \]
    so the first term on the right-hand side of (\ref{eq:shift-Phi}) is nonnegative. On the other hand, the choice of $\wh r$ maximizes the potential among the thresholds for $\wh H$ which satisfy (\ref{eq:constraint-2})--(\ref{eq:constraint-4}). Therefore, the left-hand side is nonpositive, and
    \[
        \left(\frac{3B}{4}-1\right)|\wh H_{\wh r-3}|+(B-1)|\wh H_{<\wh r-3}|\le (B-1)\sum_{i=1}^{\ell}t_i.
    \]
    Recalling that $B\ge 4$, we obtain
    \[
        |\wh H_{<\wh r-2}|=|\wh H_{\wh r-3}|+|\wh H_{<\wh r-3}|\le \frac{3}{2}\sum_{i=1}^{\ell}t_i,
    \]
    which proves the claim.

    Now suppose that $j\in\{1,2\}$ and $i\in[\ell]$, that (\ref{eq:constraint-3}) does not hold for $(j,i)$, and that $\wh H_{\wh r-j}\setminus\overline Z_{i-1}\ne\emptyset$. We show that
    \begin{equation}\label{eq:level-bound}
        |\wh H_{\wh r-j}\setminus\overline Z_{i-1}|<2\sum_{h\le i}t_h\le 2\sum_{h=1}^{\ell}t_h.
    \end{equation}
    Since (\ref{eq:constraint-4}) holds, $\wh H_{\wh r-j}\setminus\overline Z_i=\emptyset$. The failure of (\ref{eq:constraint-3}) therefore gives
    \[
        \sum_{h\le i}t_h>2^{-i}\lceil|\wh H_{\wh r-j}|\rceil_B.
    \]
    If $i=1$, then
    \[
        |\wh H_{\wh r-j}|\le \lceil|\wh H_{\wh r-j}|\rceil_B<2t_1.
    \]
    If $i>1$, condition (\ref{eq:constraint-4}) cannot have been used at the preceding step, since this would imply $\wh H_{\wh r-j}\setminus\overline Z_{i-1}=\emptyset$. Hence, (\ref{eq:constraint-3}) holds at the preceding step and
    \[
        |\wh H_{\wh r-j}\setminus\overline Z_{i-1}|+\sum_{h\le i-1}t_h
        \le 2^{-(i-1)}\lceil|\wh H_{\wh r-j}|\rceil_B
        <2\sum_{h\le i}t_h.
    \]
    This proves (\ref{eq:level-bound}).

    For each $j\in\{1,2\}$, the set $\wh H_{\wh r-j}\setminus\overline Z_{i-1}$ can be nonempty for at most one index $i$ at which (\ref{eq:constraint-3}) does not hold. Hence, if $\sum_{i=1}^{\ell}t_i>0$, the definition of the sets $\wh T_i$, (\ref{eq:large-levels}), and (\ref{eq:level-bound}) give
    \[
        \sum_{i=1}^{\ell}|\wh T_i|<\frac{3}{2}\sum_{i=1}^{\ell}t_i+4\sum_{i=1}^{\ell}t_i<6\sum_{i=1}^{\ell}t_i.
    \]
    If $\sum_{i=1}^{\ell}t_i=0$, then (\ref{eq:large-levels}) gives $\wh H_{<\wh r-2}=\emptyset$. Furthermore, if $\wh H_{\wh r-j}\setminus\overline Z_{i-1}$ were nonempty at an index where (\ref{eq:constraint-3}) does not hold, then (\ref{eq:level-bound}) would give a positive cardinality strictly less than zero. Thus, $\wh T_i=\emptyset$ for every $i\in[\ell]$.
\end{proof}

\subsection{Proof of Theorem \ref{thm:sharp-selector}}\label{subsec:proof}

We will now prove Theorem \ref{thm:coarse-selector}. Let $q=2^{17}e\,p$, and let $C_0$ be a sufficiently large absolute constant. If $q\ge1$, then $C_0pL\ge1$, and the conclusion is immediate. We may therefore assume that $q<1$.

We first bound the number of possible vectors $\mathbf{t}$.

\begin{lemma}\label{lem:possible-steps}
    Fix $\mathbf Z$ and an integer $k\ge1$. As $\mathbf t$ ranges over the size vectors satisfying $\sum_{i=1}^{\ell}t_i=k$ which arise from a minimum tower mapping to $\mathbf Z$, for each $j\in\{1,2\}$ there are at most
    \[
        3\left(\lceil\log_2(2k)\rceil+1\right)^2
    \]
    possible indices $i$ for which (\ref{eq:constraint-3}) does not hold and $\wh H_{\wh r-j}\setminus\overline Z_{i-1}\ne\emptyset$.
\end{lemma}
\begin{proof}
    For every such index $i$,
    \[
        1\le 2^{-(i-1)}\lceil|\wh H_{\wh r-j}|\rceil_B<2k.
    \]
    Indeed, (\ref{eq:constraint-4}) gives $\wh H_{\wh r-j}\setminus\overline Z_i=\emptyset$, so the failure of (\ref{eq:constraint-3}) at step $i$ gives
    \[
        2^{-i}\lceil|\wh H_{\wh r-j}|\rceil_B<\sum_{h\le i}t_h\le k,
    \]
    which gives the upper bound. The lower bound is immediate when $i=1$. When $i>1$, (\ref{eq:constraint-4}) cannot hold at step $i-1$, since $\wh H_{\wh r-j}\setminus\overline Z_{i-1}\ne\emptyset$. Hence, (\ref{eq:constraint-3}) holds at step $i-1$ and gives the lower bound.

    Let $\lceil|\wh H_{\wh r-j}|\rceil_B=B^m=2^{mL}$ and let $d=\lceil\log_2(2k)\rceil$. The previous bound implies
    \[
        mL-\log_2(2k)<i-1\le mL.
    \]
    For each fixed $m$, this interval contains at most $d+1$ integer values of $i$. If it contains a value $i\le\ell$, then
    \[
        mL<i-1+\log_2(2k)\le\ell-1+d=3L+1+d.
    \]
    Thus, the number of possible nonnegative integers $m$ is at most
    \[
        3+\left\lceil\frac{d+1}{L}\right\rceil\le3(d+1).
    \]
    The result follows.
\end{proof}

\begin{lemma}\label{lem:size-vector-count}
    For every fixed $\mathbf Z$ and integer $k\ge1$, the number of vectors $\mathbf t$ with $\sum_{i=1}^{\ell}t_i=k$ which arise from a minimum tower mapping to $\mathbf Z$ is at most $2^{10k}$.
\end{lemma}
\begin{proof}
    By Lemma \ref{lem:fragment-contain}, $t_i>0$ only if $\wh T_i\ne\emptyset$. Apart from the first step, for each $j\in\{1,2\}$ there is at most one such step, and Lemma \ref{lem:possible-steps} gives at most $3(\lceil\log_2(2k)\rceil+1)^2$ possibilities for it. Thus, there are at most
    \[
        \left(3(\lceil\log_2(2k)\rceil+1)^2+1\right)^2
    \]
    possible choices for these two steps. Once they are fixed, there are at most $\binom{k+2}{2}$ nonnegative vectors of total size $k$ supported on the resulting set of at most three steps. Since $\lceil\log_2(2k)\rceil\le k$, we have
    \[
        \left(3(\lceil\log_2(2k)\rceil+1)^2+1\right)^2\binom{k+2}{2}\le2^{10}k^6\le2^{10k}.\qedhere
    \]
\end{proof}

Let $W_1,\dots,W_\ell\sim X_q$ be independent. For each pair $(\mathbf W,H)$, we select the minimum tower of fragments using the deterministic choice fixed above. Let $\mathcal T(\mathbf W)$ be the collection of distinct selected tower tuples as $H$ varies in $\cH$. Since the fragments are pairwise disjoint, for every $\mathbf T\in\mathcal T(\mathbf W)$,
\[
    \left|\bigcup_{i=1}^{\ell}T_i\right|=\sum_{i=1}^{\ell}|T_i|.
\]

We next record the inverse multiplicity bound which follows from the two key properties in the previous subsection.
\begin{lemma}\label{lem:inverse-multiplicity}
    Fix $\mathbf Z$ and a vector $\mathbf t$ with $\sum_{i=1}^{\ell}t_i=k\ge1$. The number of selected tower tuples $\mathbf T$ which map to $(\mathbf Z,\mathbf t)$ is at most $(2^3e)^k$. Consequently, for fixed $\mathbf Z$ and $k$, the number of pairs $(\mathbf t,\mathbf T)$ arising from selected towers mapping to $\mathbf Z$ is at most $(2^{13}e)^k$.
\end{lemma}
\begin{proof}
    For the fixed pair $(\mathbf Z,\mathbf t)$, the chosen pair $(\wh r,\wh H)$ and the sets $\wh T_i$ are fixed. Lemma \ref{lem:fragment-contain} gives $T_i\subseteq\wh T_i$ for every selected tower counted above. Treating the sets $\wh T_i$ as disjoint labeled copies and using Lemma \ref{lem:fragment-contain-2}, the number of possible towers for this fixed size vector is at most
    \[
        \prod_{i=1}^{\ell}\binom{|\wh T_i|}{t_i}
        \le\binom{\sum_{i=1}^{\ell}|\wh T_i|}{k}
        \le\binom{6k}{k}
        \le(2^3e)^k.
    \]
    The second statement follows from Lemma \ref{lem:size-vector-count}, since $2^{10k}(2^3e)^k=(2^{13}e)^k$.
\end{proof}

We are now ready to state and prove the key switching lemma.
\begin{lemma}\label{lem:key}
    For every integer $a\ge0$, we have
    \[
        \E_{W_1,\dots,W_\ell\sim X_q}\left[\sum_{\substack{\mathbf T\in\mathcal T(\mathbf W)\\\sum_{i=1}^{\ell}|T_i|>a}}p^{\sum_{i=1}^{\ell}|T_i|}\right]
        \le\frac{1}{4}8^{-a}.
    \]
\end{lemma}
\begin{proof}
    We denote by $\P_q$ the product Bernoulli measure on $(2^X)^\ell$. Map each pair $(\mathbf W,\mathbf T)$ to $(\mathbf Z,\mathbf t,\mathbf T)$, where
    \[
        Z_i=W_i\cup T_i,\qquad t_i=|T_i|.
    \]
    This map is injective, since $W_i=Z_i\setminus T_i$. Furthermore, $T_i\cap W_i=\emptyset$, and hence
    \[
        \frac{\P_q(\mathbf W)}{\P_q(\mathbf Z)}=\left(\frac{1-q}{q}\right)^{\sum_{i=1}^{\ell}|T_i|}.
    \]
    Summing first over $\mathbf Z$ and using Lemma \ref{lem:inverse-multiplicity}, we obtain
    \begin{align*}
        \E_{\mathbf W}\left[\sum_{\substack{\mathbf T\in\mathcal T(\mathbf W)\\\sum_{i=1}^{\ell}|T_i|>a}}p^{\sum_{i=1}^{\ell}|T_i|}\right]
        &\le\sum_{k>a}(2^{13}e)^k\left(\frac{(1-q)p}{q}\right)^k\sum_{\mathbf Z}\P_q(\mathbf Z)\\
        &=\sum_{k>a}\left(2^{13}e\frac{(1-q)p}{q}\right)^k.
    \end{align*}
    By the choice $q=2^{17}e\,p$, we have $2^{13}e(1-q)p/q\le1/16$. Therefore, the last sum is at most
    \[
        \sum_{k>a}16^{-k}\le\frac{1}{4}8^{-a}.
    \]
\end{proof}

We will also use the following consequence of the exponential tail version of the Kahn--Kalai conjecture due to Bell \cite{Bell}. 
\begin{lemma}\label{lem:bounded-rank}
    There exists an absolute constant $C_1>0$ such that the following holds. Let $X$ be a finite set, let $0<p\le1$, let $r\ge1$, and let $\mathcal R\subseteq2^X$ be a family of nonempty sets of size at most $r$. Suppose that every cover $\mathcal C$ of $\mathcal R$ satisfies
    \[
        \sum_{C\in\mathcal C}p^{|C|}\ge\frac{1}{4}.
    \]
    Then, with probability at least $1-8^{-r}$, the set $X_{C_1pr}$ contains a member of $\mathcal R$.
\end{lemma}
\begin{proof}
    If $3p\ge1$, the conclusion is immediate after increasing $C_1$. We may therefore assume that $3p<1$.
    The family $\mathcal R$ is not $(3p)$-small. Indeed, a cover $\mathcal C$ with
    \[
        \sum_{C\in\mathcal C}(3p)^{|C|}\le\frac{1}{2}
    \]
    cannot contain the empty set and would satisfy
    \[
        \sum_{C\in\mathcal C}p^{|C|}\le\frac{1}{3}\sum_{C\in\mathcal C}(3p)^{|C|}\le\frac{1}{6}.
    \]
    The result follows from \cite[Theorem 3]{Bell}, applied with failure probability $8^{-r}$, yields the result, since 
    \[
        \log(r8^r)=3r+\log r\le4r.
    \]
\end{proof}

\begin{remk*}
    Bell's result is needed only for the exponentially small bound on the failure probability. For the version of Theorem \ref{thm:coarse-selector} in which the probability $1-2B^{-3}$ is replaced by $2/3$, it is not needed. Consequently, the same is true for the version of Theorem \ref{thm:sharp-selector} with probability $2/3$. In particular, the argument gives new proofs, independent of the proofs in \cite{PPT,PPK}, of both Talagrand's selector process conjecture and the Kahn--Kalai conjecture. 
\end{remk*}

We can now complete the proof of Theorem \ref{thm:coarse-selector}.
\begin{proof}[Proof of Theorem \ref{thm:coarse-selector}]
    Let
    \[
        \mathcal E=\left\{A\subseteq X:\min_{H\in\cH}\lambda_H(H\setminus A)\le B^{-2}\right\},
        \qquad
        Y=\bigcup_{i=1}^{\ell}W_i,
    \]
    For $Y\notin\mathcal E$, define
    \[
        \mathcal R(Y)=\{G\subseteq X\setminus Y:1\le|G|\le L,\ Y\cup G\in\mathcal E\}
    \]
    and
    \[
        \cH_0(Y)=\{H\in\cH:G\subseteq H\text{ for some }G\in\mathcal R(Y)\}.
    \]

    Suppose first that every cover $\mathcal C$ of $\cH_0(Y)$ satisfies
    \[
        \sum_{C\in\mathcal C}p^{|C|}\ge\frac{1}{4}.
    \]
    Every cover of $\mathcal R(Y)$ is then also a cover of $\cH_0(Y)$. Condition on $\mathbf W$, and let $S\sim X_{C_1pL}$ be independent. Applying Lemma \ref{lem:bounded-rank} on the ground set $X\setminus Y$, we see that $S$ contains a member of $\mathcal R(Y)$ with probability at least $1-8^{-L}$.
    Since $\mathcal E$ is increasing, whenever this happens, $Y\cup S\in\mathcal E$.

    It remains to bound the probability of the other case. In this case, fix a cover $\mathcal C$ of $\cH_0(Y)$ such that
    \[
        \sum_{C\in\mathcal C}p^{|C|}<\frac{1}{4},
    \]
    and consider $H\in\cH\setminus\cH_0(Y)$. Let $(T_1,\dots,T_\ell)$ be its selected minimum tower of fragments. We claim that
    \[
        \sum_{i=1}^{\ell}|T_i|>L.
    \]
    Indeed, otherwise put
    \[
        G=\left(\bigcup_{i=1}^{\ell}T_i\right)\setminus Y.
    \]
    We have $Y\cup G=\bigcup_{i=1}^{\ell}(W_i\cup T_i)$, so Lemma \ref{lem:feasible-capture} gives $Y\cup G\in\mathcal E$. Moreover, $G$ is nonempty since $Y\notin\mathcal E$, and $G\subseteq H$. Thus, $G\in\mathcal R(Y)$, contradicting $H\notin\cH_0(Y)$.

    The sets $\bigcup_{i=1}^{\ell} T_i$, together with $\mathcal C$, form a cover of $\cH$. Since $\cH$ is not $p$-small, the cost of the tower unions $\bigcup_{i=1}^{\ell} T_i$ is greater than $1/4$. %
Therefore,
    \[
        \sum_{\substack{\mathbf T\in\mathcal T(\mathbf W)\\\sum_{i=1}^{\ell}|T_i|>L}}p^{\sum_{i=1}^{\ell}|T_i|}>\frac{1}{4}.
    \]
    Lemma \ref{lem:key} and Markov's inequality therefore show that the probability of this case is at most $8^{-L}$. 
    Consequently,
    \[
        \P(Y\cup S\notin\mathcal E)\le2\cdot8^{-L}=2B^{-3}.
    \]

    The set $Y$ has the same distribution as $X_{q_*}$, where
    \[
        q_*=1-(1-q)^\ell\le\ell q\le2^{19}e\,pL.
    \]
    Hence, $Y\cup S$ has density at most $(2^{19}e+C_1)pL$. Taking $C_0\ge2^{19}e+C_1$ and using monotonicity proves Theorem \ref{thm:coarse-selector}, and hence Theorem \ref{thm:sharp-selector} by the reduction above.
\end{proof}

\begin{proof}[Proof of Corollary \ref{cor:sharp-selector-tail}]
    Let $m_p = \mathbb{E} Z_{\Lambda}(X_p)$. We may assume that $0<m_p<\infty$. Let
    \[
        L=\left\lceil\log(4/\delta)\right\rceil,
        \qquad B=2^L,
        \qquad q=\frac{p}{C_0L}.
    \]
    Suppose, for the sake of contradiction, that
    \[
        \mathcal A=\{A\subseteq X:Z_\Lambda(A)>(1+\delta)m_p\}
    \]
    is not $q$-small. For every $A\in\mathcal A$, choose $\lambda^A\in\Lambda$ such that $\lambda^A(A)>(1+\delta)m_p$, and put
    \[
        \nu_A(x)=\frac{\lambda^A(x)\mathbb I(x\in A)}{\lambda^A(A)}.
    \]
    For $\nu_A(x)>0$, let $\lambda'_A(x)=B^{\lfloor\log_B\nu_A(x)\rfloor}$, and let $\lambda'_A(x)=0$ otherwise. Then
    \[
        \lambda'_A(x)\le\nu_A(x)\le B\lambda'_A(x),
        \qquad
        \frac{1}{B}\le\lambda'_A(A)\le1.
    \]

    Applying Theorem \ref{thm:coarse-selector} to $\mathcal A$ at density $q$, we have, with probability at least $1-2B^{-3}$, some $A\in\mathcal A$ such that
    \[
        \lambda'_A(A\setminus X_p)\le B^{-2}.
    \]
    In this case,
    \[
        \nu_A(A\setminus X_p)\le B\lambda'_A(A\setminus X_p)\le B^{-1}.
    \]
    For this $A$,
    \[
        Z_\Lambda(X_p)\ge\lambda^A(X_p)\ge\lambda^A(A\cap X_p)>
        (1-B^{-1})(1+\delta)m_p.
    \]
    Taking expectations gives
    \[
        m_p>(1-2B^{-3})(1-B^{-1})(1+\delta)m_p.
    \]
    On the other hand, $B^{-1}\le\delta/4$ and $2B^{-3}\le\delta/4$, while
    \[
        (1-\delta/4)^2(1+\delta)
        =1+\frac{\delta}{16}(8-7\delta+\delta^2)>1
    \]
    for $0<\delta\le1$, a contradiction. Thus, $\mathcal A$ is $q$-small, from which the desired result follows. %
\end{proof}

\section{Finding small integral covers from fractional covers}\label{sec:fractional-covers}

In this section, we give the proof of Theorem \ref{thm:main}, relying on Theorem \ref{thm:sharp-selector}. In particular, given a fractional cover $w$, we construct appropriate weight vectors $\lambda_H$ on $H$. If $\cH$ is not $q$-small, then Theorem \ref{thm:sharp-selector} shows that, with reasonably high probability, $\max_{H\in \cH} \lambda_{H}(X_{q'})$ is very close to $1$, for $q'$ not much larger than $q$. Relying on this conclusion, we then show that the sum of $w(W)$ over $W$ contained in $X_{q'}$ must be large. However, if $q' < p$, then this contradicts the assumption that $w$ is fractionally small, i.e. $\sum_{W\in 2^X}w(W)p^{|W|}\le1/2$.

\begin{proof} [Proof of Theorem \ref{thm:main}]
We may first assume that $w(\emptyset)=0$. Indeed, if $a=w(\emptyset)>0$, define $w'(\emptyset)=0$ and
\[
    w'(W)=\min\left\{\frac{w(W)}{1-a},1\right\}
    \qquad\textrm{for every nonempty }W\subseteq X.
\]
Since $a\le1/2$, the function $w'$ is a fractional cover of $\cH$, and
\[
    \sum_{W\in2^X}w'(W)p^{|W|}
    \le\frac{\sum_{W\in2^X}w(W)p^{|W|}-a}{1-a}
    \le\frac{1/2-a}{1-a}\le\frac12.
\]
Indeed, for each $H\in\cH$, the sum of $w(W)/(1-a)$ over the nonempty sets $W\subseteq H$ is at least $1$, and truncating a term can only occur when that term becomes $1$. We can therefore replace $w$ by $w'$.

Given the fractional cover $w$, we now construct weight vectors $\lambda_H$ for $H\in \cH$ as follows. Denote $\cW$ the support of the fractional cover $w$. For each $H\in \cH$, define 
\[
    \lambda_H(x) = \frac{\sum_{W \in \cW: x\in W, W\subseteq H} \frac{w(W)}{|W|}}{\sum_{y\in H}\sum_{W \in \cW: y\in W, W\subseteq H} \frac{w(W)}{|W|}}. 
\]
In particular, $\sum_{x\in H}\lambda_H(x) = 1$. Also note that 
\[
    \sum_{y\in H} \,\,\,\,\, \sum_{W \in \cW: y\in W, W\subseteq H} \frac{w(W)}{|W|} = \sum_{W\in \cW: W\subseteq H} w(W).
\]

Let $s=\lceil\log(2t)\rceil$. Assume for the sake of contradiction that $\cH$ is not $q$-small for
\[
    q=\frac{cp}{\log(2t)}.
\]
By Theorem \ref{thm:sharp-selector}, with probability at least $2/3$, for $q'=Csq$ and an absolute constant $C>0$,
    \[
        \max_{H\in \cH} \sum_{x\in X_{q'}\cap H} \lambda_H(x) \ge 1-2^{-s}\ge 1 - \frac{1}{2t}.
    \]
Since $s\le2\log(2t)$, we have $q'\le2Ccp$. We choose $c\le\min\{1,(12C)^{-1}\}$, so that $q\le1$ and $q'\le p/6$.
On the other hand, observe that 
\begin{align*}
    \sum_{x\in X_{q'}\cap H} \lambda_H(x) &= \frac{\sum_{W\in \cW:W\subseteq H}w(W) - \sum_{W\in \cW:W\subseteq H, W\not \subseteq X_{q'}}w(W)\frac{|W\setminus X_{q'}|}{|W|}}{\sum_{W\in \cW:W\subseteq H}w(W)}\\
    &= 1 - \frac{\sum_{W\in \cW:W\subseteq H, W\not \subseteq X_{q'}}w(W)\frac{|W\setminus X_{q'}|}{|W|}}{\sum_{W\in \cW:W\subseteq H\cap X_{q'}}w(W) + \sum_{W\in \cW:W\subseteq H, W\not\subseteq X_{q'}}w(W)}\\
    &\le 1 - \frac{\frac{1}{t}\sum_{W\in \cW:W\subseteq H, W\not\subseteq X_{q'}}w(W)}{\sum_{W\in \cW:W\subseteq H\cap X_{q'}}w(W) + \sum_{W\in \cW:W\subseteq H, W\not\subseteq X_{q'}}w(W)}\\
    &= 1 - \frac{1}{t} + \frac{\sum_{W\in \cW:W\subseteq H\cap X_{q'}}w(W)}{t\sum_{W\in \cW:W\subseteq H}w(W)}.
\end{align*}
Here, in the first inequality, we use that $|W|\le t$ for all $W\in \cW$. 

We hence obtain, with probability at least $2/3$,
\begin{align*}
    \frac{1}{2t} &\le -1 + \frac{1}{t} + \max_{H\in \cH} \sum_{x\in X_{q'}\cap H} \lambda_H(x) \\
    &\le \max_{H\in \cH}\frac{\sum_{W\in \cW:W\subseteq H\cap X_{q'}}w(W)}{t\sum_{W\in \cW:W\subseteq H}w(W)}\\
    &\le \max_{H\in \cH}\frac{\sum_{W\in \cW:W\subseteq H\cap X_{q'}}w(W)}{t}.
\end{align*}
Here we use the assumption that $w$ is a fractional cover of $\cH$ and hence $\sum_{W\in \cW:W\subseteq H}w(W) \ge 1$. Therefore, with probability at least $2/3$
\begin{align*}
    \max_{H\in \cH}\sum_{W\in \cW:W\subseteq H\cap X_{q'}}w(W) \ge \frac{1}{2},
\end{align*}
which implies
\begin{align*}
    \E\lt[\max_{H\in \cH}\sum_{W\in \cW:W\subseteq H\cap X_{q'}}w(W)\rt] \ge \frac{1}{3},
\end{align*}

Note that 
\[
    \sum_{W\in \cW:W\subseteq X_{q'}} w(W) \ge \max_{H\in \cH}\sum_{W\in \cW:W\subseteq H\cap X_{q'}}w(W).
\]
Thus,
\begin{align*}
    \E\left[\sum_{W\in \cW:W\subseteq X_{q'}} w(W)\right] \ge \frac{1}{3}.
\end{align*}

On the other hand, by linearity of expectation, 
\begin{align*}
    \E\left[\sum_{W\in \cW:W\subseteq X_{q'}} w(W)\right] = \sum_{W\in \cW} (q')^{|W|}w(W). 
\end{align*}
Since $q'\le p/6$ and every set in $\cW$ is nonempty,
\begin{align*}
    \E\left[\sum_{W\in \cW:W\subseteq X_{q'}} w(W)\right]
    &\le \sum_{W\in \cW} (p/6)^{|W|} w(W)\\
    &\le \frac16\sum_{W\in \cW}p^{|W|}w(W)\le\frac{1}{12}.
\end{align*}
This is a contradiction, and hence we conclude that $\cH$ must be $q$-small, as desired.  
\end{proof}

Combining Theorem \ref{thm:main} with the sampling argument of Fischer and Person \cite[Lemmas 16 and 17]{FP}, which we include below for completeness, we obtain Corollary \ref{corl:loglog-gap}. 

\begin{proof}[Proof of Corollary \ref{corl:loglog-gap}]
    Let $n=|X|$ and $t=\lceil\log(8n)\rceil$. As in the proof of Theorem \ref{thm:main}, we may assume that $w(\emptyset)=0$. Put
    \[
        \cH_0=\left\{H\in\cH:
        \sum_{\substack{W\subseteq H\\1\le |W|\le t}}w(W)\ge\frac12\right\},
        \qquad
        \cH_1=\cH\setminus\cH_0.
    \]
    For $1\le |W|\le t$, let $w_0(W)=\min\{2w(W),1\}$, and let $w_0(W)=0$ otherwise. Then $w_0$ is a fractional cover of $\cH_0$ and
    \[
        \sum_{W\in2^X}w_0(W)(p/2)^{|W|}
        \le\sum_{W\in2^X}w(W)p^{|W|}\le\frac12.
    \]
    Let $c_0\le1$ be the constant in Theorem \ref{thm:main}, and set
    \[
        q=\frac{c_0p}{4\log(2t)}.
    \]
    Applying Theorem \ref{thm:main} to $w_0$ at density $p/2$, we obtain a cover of $\cH_0$ having cost at most $1/2$ at density $c_0p/(2\log(2t))$. This cover cannot contain the empty set. Decreasing the density by a factor of $2$, we therefore obtain a cover $\cG_0$ of $\cH_0$ such that
    \begin{equation}\label{eq:small-part-loglog}
        \sum_{G\in\cG_0}q^{|G|}\le\frac14.
    \end{equation}

    Independently include every $W\subseteq X$ with $|W|>t$ in $\cG_1$ with probability $\min\{4(n+1)w(W),1\}$. For every $H\in\cH_1$,
    \[
        \sum_{\substack{W\subseteq H\\|W|>t}}w(W)>\frac12.
    \]
    If one of these sets is included with probability $1$, then $H$ is covered surely. Otherwise,
    \[
        \P(H\notin\langle\cG_1\rangle)
        \le\exp\left(-4(n+1)\sum_{\substack{W\subseteq H\\|W|>t}}w(W)\right)
        <e^{-2(n+1)}.
    \]
    Hence,
    \[
        \P(\cH_1\nsubseteq\langle\cG_1\rangle)
        \le2^ne^{-2(n+1)}<\frac14.
    \]
    Since $q\le p/4$ and $t\ge\log(8n)$,
    \begin{align*}
        \E\left[\sum_{G\in\cG_1}q^{|G|}\right]
        &\le4(n+1)\sum_{|W|>t}w(W)(p/4)^{|W|}\\
        &\le2(n+1)4^{-(t+1)}\le\frac1{64}.
    \end{align*}
    By Markov's inequality, with probability at least $15/16$,
    \[
        \sum_{G\in\cG_1}q^{|G|}\le\frac14.
    \]
    The event that $\cG_1$ covers $\cH_1$ has probability greater than $3/4$, while the last event has probability at least $15/16$. Their intersection is nonempty, so we may fix a choice of $\cG_1$ which covers $\cH_1$ and satisfies the last inequality. Together with (\ref{eq:small-part-loglog}), this shows that $\cH$ is $q$-small.
\end{proof}

\section{Sudakov minoration for selector processes}\label{sec:sudakov}

In this section, we prove Theorem \ref{thm:sudakov} and its corollaries.

\begin{proof}[Proof of Theorem \ref{thm:sudakov}]
    For $T\subseteq X$, let
    \[
        \Lambda_T=\{\lambda\in\Lambda_0:T\subseteq C_\lambda\},
        \qquad
        d(T)=|\Lambda_T|.
    \]
    Choose $T\subseteq X$ maximizing
    \begin{equation}\label{eq:sudakov-maximizer}
        d(T)p^{-|T|},
    \end{equation}
    and, among all maximizers, choose one for which $|T|$ is minimum. Since
    \[
        d(\emptyset)=|\Lambda_0|\ge p^{-n},
    \]
    the maximum is positive, and hence $T\subseteq C_\lambda$ for some $\lambda\in\Lambda_0$. In particular, $|T|\le n$.

    We first observe that
    \begin{equation}\label{eq:sudakov-two-vectors}
        d(T)\ge2.
    \end{equation}
    Indeed, if $d(T)=1$, then
    \[
        d(T)p^{-|T|}\le p^{-n}\le d(\emptyset).
    \]
    Thus, the empty set is also a maximizer in (\ref{eq:sudakov-maximizer}). By the choice of $T$, we would have $T=\emptyset$, which contradicts $d(T)=1$, since $d(\emptyset)\ge p^{-n}>1$.

    By the maximality of $T$, for every $S\subseteq X\setminus T$,
    \begin{equation}\label{eq:sudakov-density}
        d(T\cup S)\le p^{|S|}d(T).
    \end{equation}
    Consider the family
    \[
        \cH_T=\{C_\lambda\setminus T:\lambda\in\Lambda_T\}.
    \]
    Its members are distinct. Indeed, if $C_\lambda\setminus T=C_{\lambda'}\setminus T$ for distinct $\lambda,\lambda'\in\Lambda_T$, then $C_\lambda=C_{\lambda'}$, contradicting (\ref{eq:sudakov-private-mass}).

    We claim that $\cH_T$ is not $p$-small. Let $\cG$ be a cover of $\cH_T$, after discarding from $\cG$ any set which is not contained in a member of $\cH_T$. Counting the members of $\Lambda_T$ according to a set of the cover contained in $C_\lambda\setminus T$, and then using (\ref{eq:sudakov-density}), we obtain
    \[
        d(T)
        \le\sum_{G\in\cG}d(T\cup G)
        \le d(T)\sum_{G\in\cG}p^{|G|}.
    \]
    Thus, every cover $\cG$ of $\cH_T$ satisfies
    \[
        \sum_{G\in\cG}p^{|G|}\ge1,
    \]
    proving the claim.

    For $H=C_\lambda\setminus T\in\cH_T$, choose, using (\ref{eq:sudakov-two-vectors}), some $\lambda'\in\Lambda_T$ distinct from $\lambda$. Since $T\subseteq C_{\lambda'}$, we have
    \[
        \lambda(H)
        \ge\lambda(C_\lambda\setminus C_{\lambda'})
        \ge r.
    \]
    Define a weight vector $\nu_H:X\to[0,1]$ by
    \[
        \nu_H(x)=
        \begin{cases}
            \lambda(x)/\lambda(H), & x\in H,\\
            0, & x\notin H.
        \end{cases}
    \]
    Then $\nu_H$ is supported on $H$ and $\nu_H(H)=1$.

    Applying Theorem \ref{thm:sharp-selector} to the family $\cH_T$ and the weight vectors $\nu_H$, we have, with probability at least $1-2^{-s-1}$, some $H=C_\lambda\setminus T\in\cH_T$ such that
    \[
        \nu_H(X_{Csp})\ge1-2^{-s}.
    \]
    Since $\lambda(x)=\lambda(H)\nu_H(x)\ge r\nu_H(x)$ for every $x\in H$, it follows that
    \[
        Z_\Lambda(X_{Csp})
        \ge\lambda(X_{Csp})
        \ge r\nu_H(X_{Csp})
        \ge(1-2^{-s})r.
    \]
\end{proof}

In particular, let $C_0\ge1$ be the absolute constant in Theorem \ref{thm:sudakov}, and choose
\[
    s=\lceil\log(1/\eps)\rceil,
    \qquad
    q=\frac{p}{C_0s}.
\]
Since $s\le\log(2/\eps)$, the hypothesis in (\ref{eq:sudakov-fixed-density}), with $C=C_0$, gives $|\Lambda_0|\ge q^{-n}$. Theorem \ref{thm:sudakov}, applied at density $q$, gives
\[
    \P\bigl(Z_\Lambda(X_p)\ge(1-\eps)r\bigr)
    \ge1-2^{-s-1}
    \ge1-\frac{\eps}{2}.
\]

\begin{proof}[Proof of Corollary \ref{corl:sudakov-expectation}]
    Let $C\ge 1$ be the absolute constant in Theorem \ref{thm:sudakov}, and let $K=\lceil C\rceil$. Applying Theorem \ref{thm:sudakov} with $s=1$ gives
    \begin{equation}\label{eq:sudakov-expectation-large-density}
        \E Z_\Lambda(X_{\min\{Cp,1\}})\ge\frac{3r}{8}.
    \end{equation}
    Indeed, the corresponding random variable is at least $r/2$ with probability at least $3/4$.

    Let $Y\sim X_{\min\{Cp,1\}}$. Independently assign each element of $Y$ one of $K$ colors uniformly, and let $Y_i$ be the set of elements receiving color $i$. Each $Y_i$ has the same law as $X_\rho$, where
    \[
        \rho=\frac{\min\{Cp,1\}}{K}\le p.
    \]
    Since the weight vectors are nonnegative,
    \[
        Z_\Lambda(Y)
        \le\sum_{i=1}^K Z_\Lambda(Y_i).
    \]
    Taking expectations, using monotonicity, and then (\ref{eq:sudakov-expectation-large-density}), we get
    \[
        \frac{3r}{8}
        \le\E Z_\Lambda(Y)
        \le K\E Z_\Lambda(X_\rho)
        \le K\E Z_\Lambda(X_p).
    \]
    This proves the result with $c=3/(8K)$.
\end{proof}

\begin{proof}[Proof of Corollary \ref{corl:talagrand-sudakov}]
    If $p=1$, the result is immediate. Suppose that $p<1$, and take
    \[
        \Lambda=\Lambda_0=\{\mathbb I_F:F\in\cF\},
        \qquad
        C_{\mathbb I_F}=F.
    \]
    For distinct $F,F'\in\cF$,
    \[
        \mathbb I_F(C_{\mathbb I_F}\setminus C_{\mathbb I_{F'}})
        =|F\setminus F'|
        =n-|F\cap F'|
        \ge\frac n2.
    \]
    Corollary \ref{corl:sudakov-expectation}, applied with $r=n/2$, gives
    \[
        \E\max_{F\in\cF}|F\cap X_p|
        =\E Z_\Lambda(X_p)
        \ge cn,
    \]
    after changing the absolute constant $c$.
\end{proof}

\section*{Acknowledgement}
The author would like to thank Jan Vondr\'ak for stimulating discussions around the topic of the paper.

\end{document}